\documentclass[12pt]{article}
\usepackage{amssymb}
\usepackage{amsmath}
\usepackage{graphicx}
\usepackage{color}

\newcommand{\powerset}[1]{\ensuremath{\mathcal{P}(#1)}}
\newcommand{\lattice}[1]{\ensuremath{\mathcal{L}(#1)}}
\newcommand{\polyVenn}[1]{$#1$-polyVenn}
\newcommand{\polyVenns}[1]{$#1$-polyVenns}

\begin{document}
\begin{center}
\Large
Minimum Area Venn Diagrams Whose Curves are Polyominoes
\end{center}

\begin{flushright}
Stirling Chow and Frank Ruskey \\
Department of Computer Science \\
University of Victoria \\
Victoria, B.C., Canada  V8W 3P6 \\
\verb+schow@cs.uvic.ca+
\end{flushright}

While working at the Berlin Academy, the renowned Swiss
  mathematician Leonard Euler was asked to tutor Frederick the
  Great's niece, the Princess of Anhalt-Dessau, in all matters of
  natural science and philosophy \cite{James}.
Euler's tutelage of the princess continued from $1760$ to $1762$
  and culminated in the publishing of the popular and
  widely-translated ``Letters to a German Princess'' \cite{Euler}.
In the letters, Euler eloquently wrote about diverse topics
  ranging from why the sky was blue to free will and determinism.

In his lesson on categorical propositions and syllogisms, Euler
  used diagrams comprised of overlapping circles; these diagrams
  became known as Eulerian circles, or simply Euler diagrams.
In an Euler diagram, a proposition's classes are represented as
  circles whose overlap depends on the relationship established by
  the proposition.
For example, the propositions

\begin{figure}
\centerline{\includegraphics[scale=0.9]{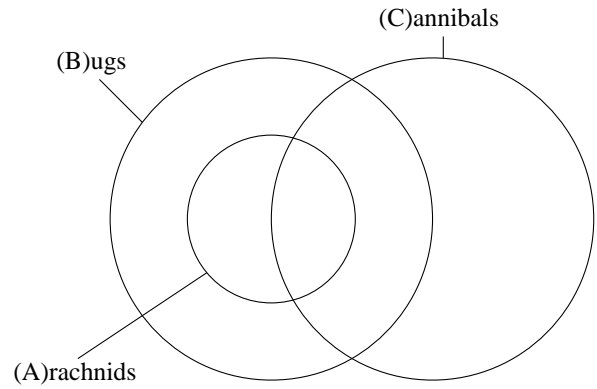}} \caption[ ]{An
example of an Euler diagram.} \label{fig:euler}
\end{figure}

\begin{center}
\begin{tabular}{l} All arachnids are bugs \\
Some bugs are cannibals
\end{tabular}
\end{center}

\noindent can be represented by Fig. \ref{fig:euler}.

\begin{figure}
\centering \scalebox{1.0}{\input{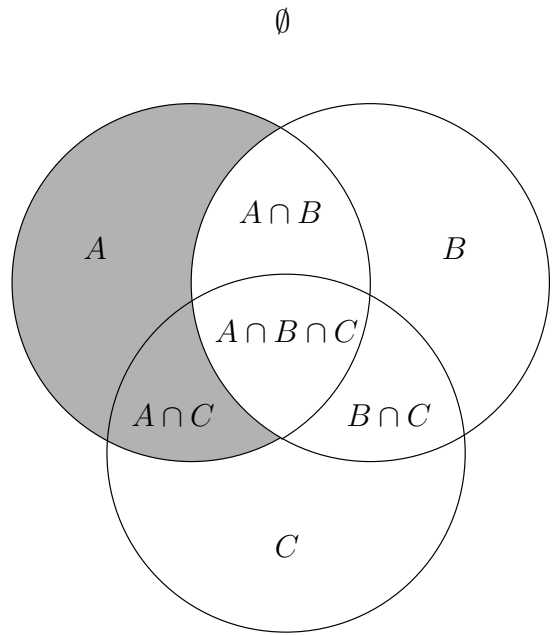}} \caption[ ]{A
Venn diagram that represents the Euler diagram in Fig.
\ref{fig:euler} by shading the missing regions.} \label{fig:venn3}
\end{figure}

In $1880$, a Cambridge priest and mathematician named John Venn
  published a paper studying special instances of Euler diagrams in
  which the classes overlap in all possible ways \cite{Venn};
  although originally applied to logic reasoning, these Venn
  diagrams are now commonly used to teach students about set theory.
For example, the Venn diagram in Fig. \ref{fig:venn3} shows
  all the ways in which three sets can intersect.
The primary difference between Venn and Euler diagrams is how they
  represent empty sets (e.g., the set of arachnids which are
  \emph{not} bugs in the example of Fig. \ref{fig:euler}).
In an Euler diagram, regions representing empty sets are
  omitted, while in Venn diagrams they are included but
  denoted by shading.

Informally, an $n$-Venn diagram is a set of $n$ simple, closed
  curves that subdivide the plane into $2^n$ connected   regions with
  each region uniquely mapping to a subset of the $n$ curves
  consisting of those curves which enclose it.
The regions are usually referred to by their enclosing curves.
For example, the $3$-Venn diagram in Fig. \ref{fig:venn3} with
  curves $\{A, B, C\}$ has regions $\{\emptyset, A, B, C, AB, AC,
  BC, ABC\}$.
If an $n$-Venn diagram's curves are equivalent to each other
  modulo translations, rotations, and reflections, then the diagram
  is referred to as a \emph{congruent} $n$-Venn diagram.

\begin{figure}
\centerline{\includegraphics[scale=0.65]{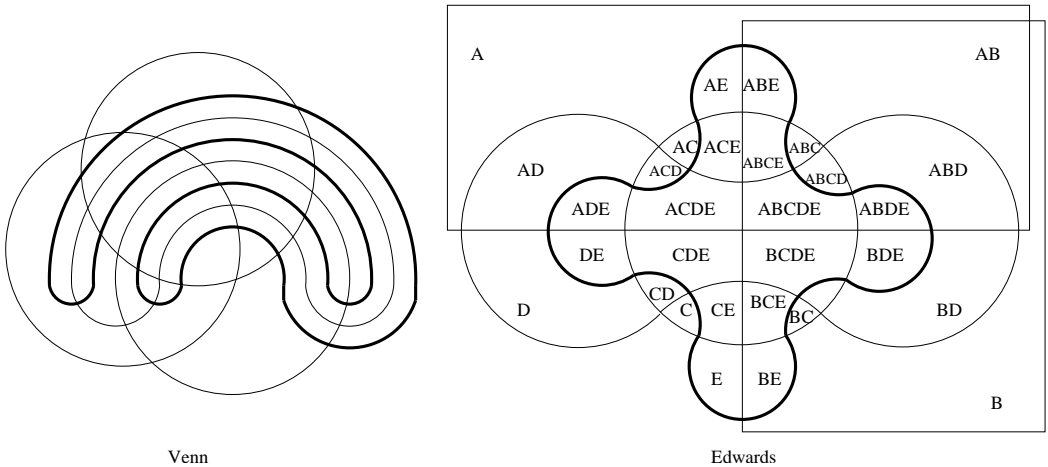}}
\caption[ ]{A $5$-Venn diagram drawn iteratively using Venn's and
Edwards' algorithms; the last curve drawn is highlighted.}
\label{fig:venn5_venn_edwards}
\end{figure}

In recent years, there has been renewed interest in studying the
  combinatorial and geometric properties of Venn diagrams
  \cite{Edwards2, Ruskey}.
Of paramount importance is how to draw a Venn diagram for a given
  number of sets.
John Venn proposed an iterative algorithm in his original Venn
  paper \cite{Venn}; unfortunately, the resulting drawings lacked an
  aesthetic appeal.
In $1989$, Anthony Edwards developed an elegant method for
  drawing $n$-Venn diagrams that produced highly symmetric drawings
  \cite{Edwards1}.
Figure \ref{fig:venn5_venn_edwards} shows a comparison of $5$-Venn
  diagrams drawn using Venn's and Edwards' algorithms.

\begin{figure}
\centerline{\includegraphics[scale=0.80]{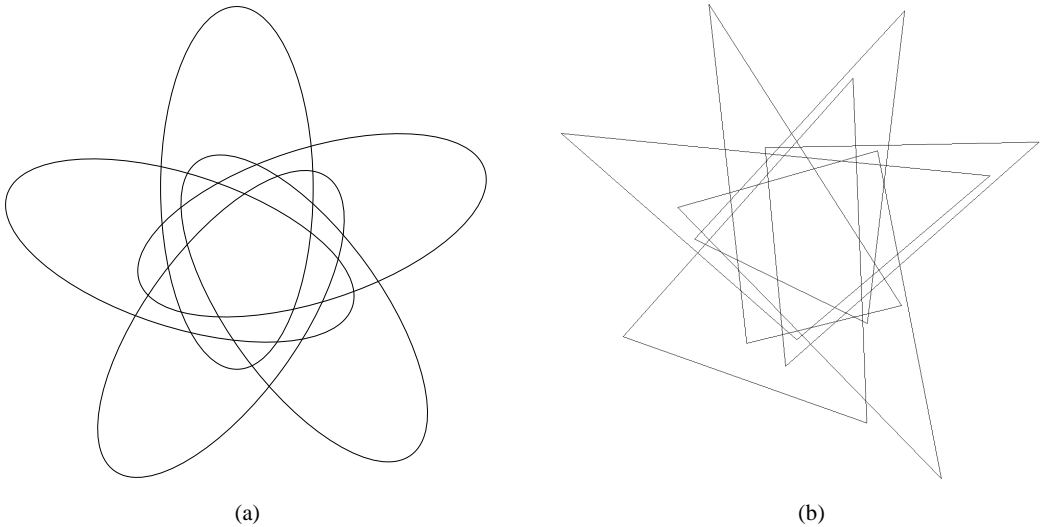}} \caption[
]{(a) A symmetric $5$-Venn diagram using ellipses and (b) a
$6$-Venn diagram using triangles.} \label{fig:venn_shapes}
\end{figure}

An interesting problem popularized by Gr\"{u}nbaum
  \cite{Grunbaum0, Grunbaum1, Grunbaum2, Grunbaum3} is to consider which Venn
  diagrams can be drawn using specific shapes.
Figure \ref{fig:venn3} shows a $3$-Venn diagram comprised of
  circles; a natural question to ask is if such a diagram exists for
  four sets.
It turns out the answer is no, and this can be proved easily using
  Euler's formula ($V = E - F + 2$) and noting that such a diagram
  must have $2^4=16$ faces and that two circles can intersect at at
  most two points \cite{Ruskey}.
Figure \ref{fig:venn_shapes} shows examples of Venn diagrams
  drawn using ellipses \cite{Grunbaum0} and triangles \cite{Carroll}.
The diagram in Fig. \ref{fig:venn_shapes}(a) is special because it
  is an example of a \emph{symmetric} Venn diagram;
  that is, a diagram with $n$-fold rotational symmetry and
  (necessarily) congruent curves.
Symmetric Venn diagrams exist if and only if $n$ is prime
  \cite{Griggs}.

\begin{figure}
\centerline{\includegraphics[scale=1.0]{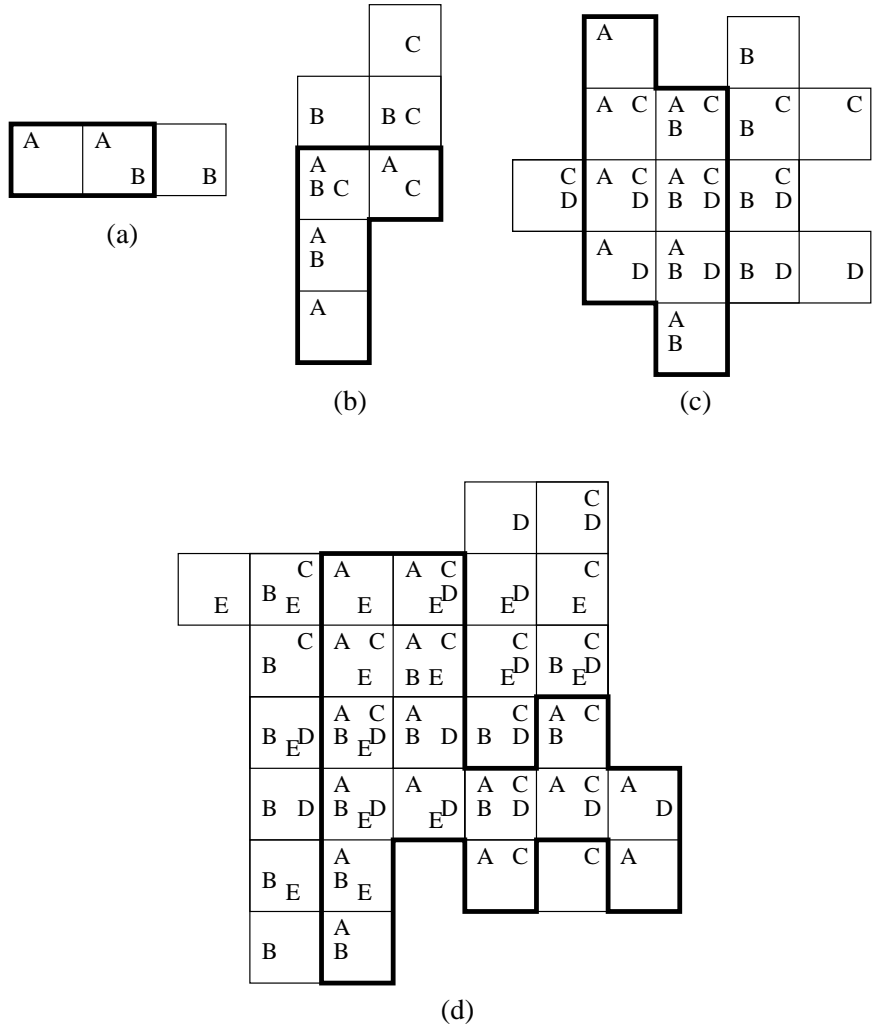}}
\caption[ ]{(a)--(c) Mark Thompson's congruent \polyVenns{n}\ for
$n=2,3,4$ and (d) the authors' congruent \polyVenn{5}; in each
case, curve $A$ is highlighted.} \label{fig:venn_polyominoes}
\end{figure}

On his ``Math Recreations'' web site \cite{Thompson}, Mark
  Thompson proposed the novel problem of finding
  \emph{Venn polyominoes}
  (from now on referred to as \emph{\polyVenns{n}});
  these are Venn diagrams whose curves are the outlines of
  polyominoes.
Polyominoes, or $n$-ominoes, are a generalization of dominoes
  ($2$-ominoes) whereby shapes are formed by gluing together $n$
  unit squares.
One can also think of a polyomino as being the result of cutting a
  shape from a piece of graph paper where the cuts are made along
  the lines.
Thompson found examples of congruent \polyVenns{n}\ for
  $n=2,3,4$, and using a computer search, we found
  a congruent \polyVenn{5}
  (see Fig. \ref{fig:venn_polyominoes}).

In the remainder of this paper, our focus is on
  minimizing the total area
  of the drawing (relative to a scaling factor).
We present some examples that minimize area according
  to various additional constraints.
At present, these examples do not generalize and
  so we develop an algorithm that comes
  close to minimizing the area.
The algorithm is simple and utilizes symmetric chain
  decompositions of the Boolean lattice.
We also provide asymptotic results that relate the area required
  by the algorithm's diagrams to the theoretical minimum area.
We conclude by presenting some open problems related to Venn
  polyominoes and other shape-constrained Venn diagrams.

\subsection*{Polyominoes}
\begin{figure}
\centerline{\includegraphics[scale=1.0]{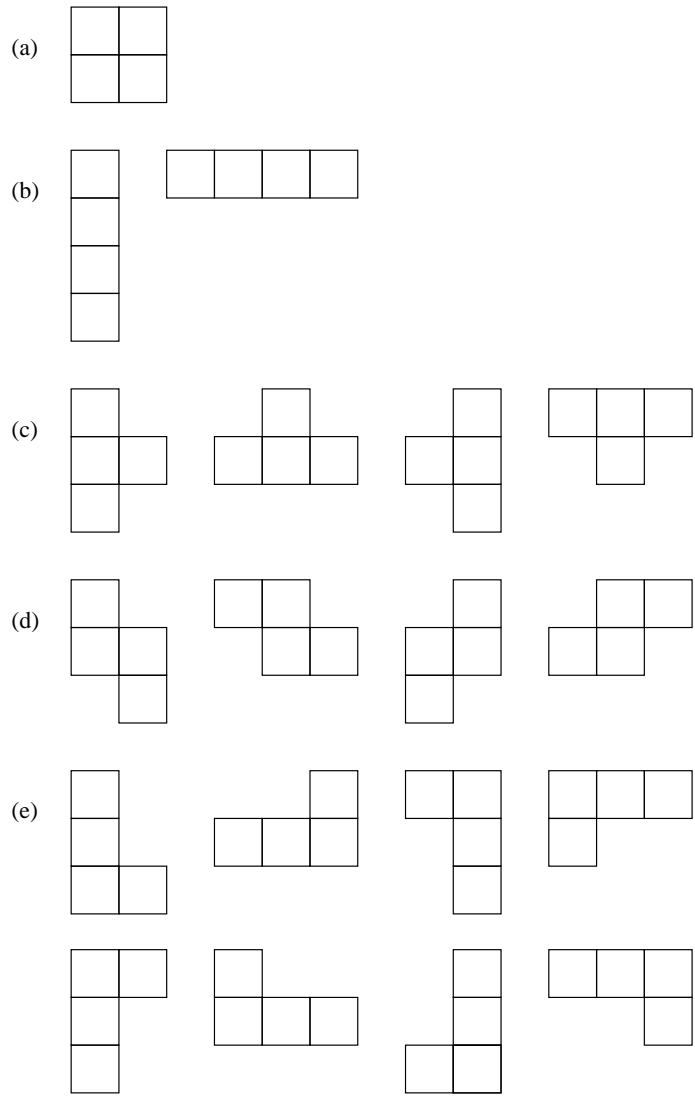}} \caption{All
possible $4$-ominoes (tetrominoes)} \label{fig:tetrominoes}
\end{figure}

A \emph{polyomino} is an edge-connected set of unit squares,
  called \emph{cells}, embedded in the integer lattice.
Two cells are adjacent if, and only if, they share a common edge.
Edge-connected means that every pair of cells is connected by
  a path through adjacent cells.
Polyominoes are often classified by area and referred to as
  $n$-ominoes when they contain $n$ cells.
For example, the games of dominoes and Tetris are played with
  $2$-ominoes and $4$-ominoes (tetrominoes), respectively (see
  Fig. \ref{fig:tetrominoes}).

Polyominoes have been extensively studied and have a wide-range of
  applications in mathematics and the physical sciences
  \cite{Golomb, Klarner}.
The problem of counting $n$-ominoes has garnered considerable
  interest \cite{Jensen,Mertens,Redelmeier}, and although
  counts up to $47$-ominoes are known (see sequence $A001168$
  \cite{Sloane}), the problem of finding a formula
  for the number remains open.

Several other subclasses of polyominoes have been defined.
\emph{Free} polyominoes treat polyominoes that are translations,
  rotations, or reflections of each other to be equivalent whereas
  \emph{fixed} polyominoes only consider translations as being
  equivalent.
For example, Fig. \ref{fig:tetrominoes} shows the $19$
  equivalence classes of fixed tetrominoes and $5$ equivalence
  classes (a,b,c,d, and e), of free tetrominoes.

\begin{figure}
\centerline{\includegraphics[scale=1.0]{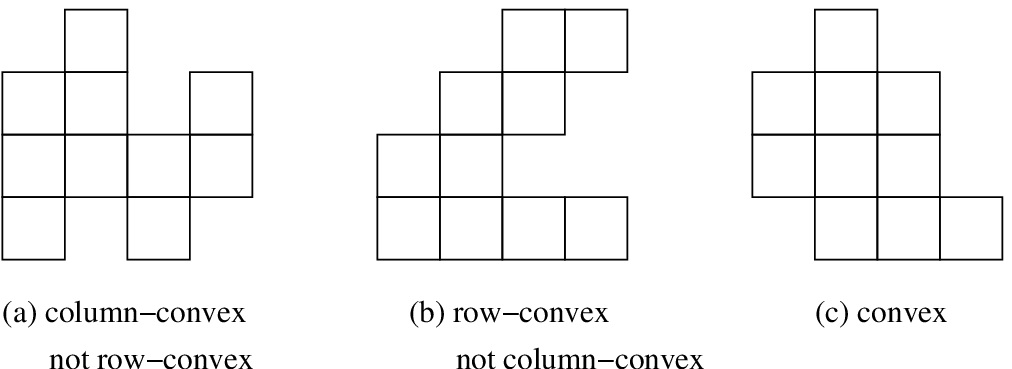}}
\caption{$10$-ominoes that exhibit different convexivity
properties.} \label{fig:convexivity}
\end{figure}

If every column (row) of a polyomino is a contiguous strip of
  cells then the polyomino is called \emph{column-convex}
  (\emph{row-convex}).
A \emph{convex} polyomino is one that is both column and row
  convex (see Fig. \ref{fig:convexivity}).
No closed-form formula is known for the number, $a(n)$, of
  fixed column-convex $n$-ominoes; however, P\'{o}lya \cite{Polya}
  derived the recurrence relation
  $a(n)=5a(n-1)-7a(n-2)+4a(n-3)$ with $a(1)=1$, $a(2)=2$, $a(3)=6$,
  and $a(4)=19$.
This recurrence relation has the rational generating function
  $g(x) = \frac{x(1-x)^3}{1-5x+7x^2-4x^3}$ (see sequence $A001169$
  \cite{Sloane}).

\subsection*{Minimum Area \polyVenns{n}}

An \polyVenn{n}\ is a Venn diagram comprised of $n$ curves,
  each of which is the perimeter of some polyomino.
In particular, each polyomino must be free of holes in order for
  the perimeter to be a simple, closed curve, and when placed on top
  of another polyomino, may not partially cover any of the bottom
  polyomino's cells (i.e., the corners of the curves must have
  unit coordinates).

Referring to the examples in Fig. \ref{fig:venn_polyominoes}, we
  see that an \polyVenn{n}\ can be drawn by tracing the curves on
  the lines of a piece of graph paper;
  in the (combinatorial) graph
  drawing community, this is referred to as an
  \emph{orthogonal grid drawing} \cite{Battista}.
In fact, any orthogonal grid drawing of a Venn diagram will
  produce curves that are the perimeters of polyominoes.
Since each bounded region must contain at least one cell and
 there is exactly one unbounded region, the minimum area for such a
 diagram is $2^n - 1$ cells.
In addition, since each curve encloses $2^{n-1}$ regions, it must
  be the perimeter of at least a $2^{n-1}$-omino.
This leads us to the following definition of a minimum area
  \polyVenn{n}:

\vspace{12pt}
\textsc{Definition} \linebreak
\indent A \textit{minimum area \polyVenn{n}} is an orthogonal
  unit-grid drawing of a Venn diagram with area $2^n - 1$.
\vspace{12pt}

\begin{figure}
\centerline{\includegraphics[scale=0.8]{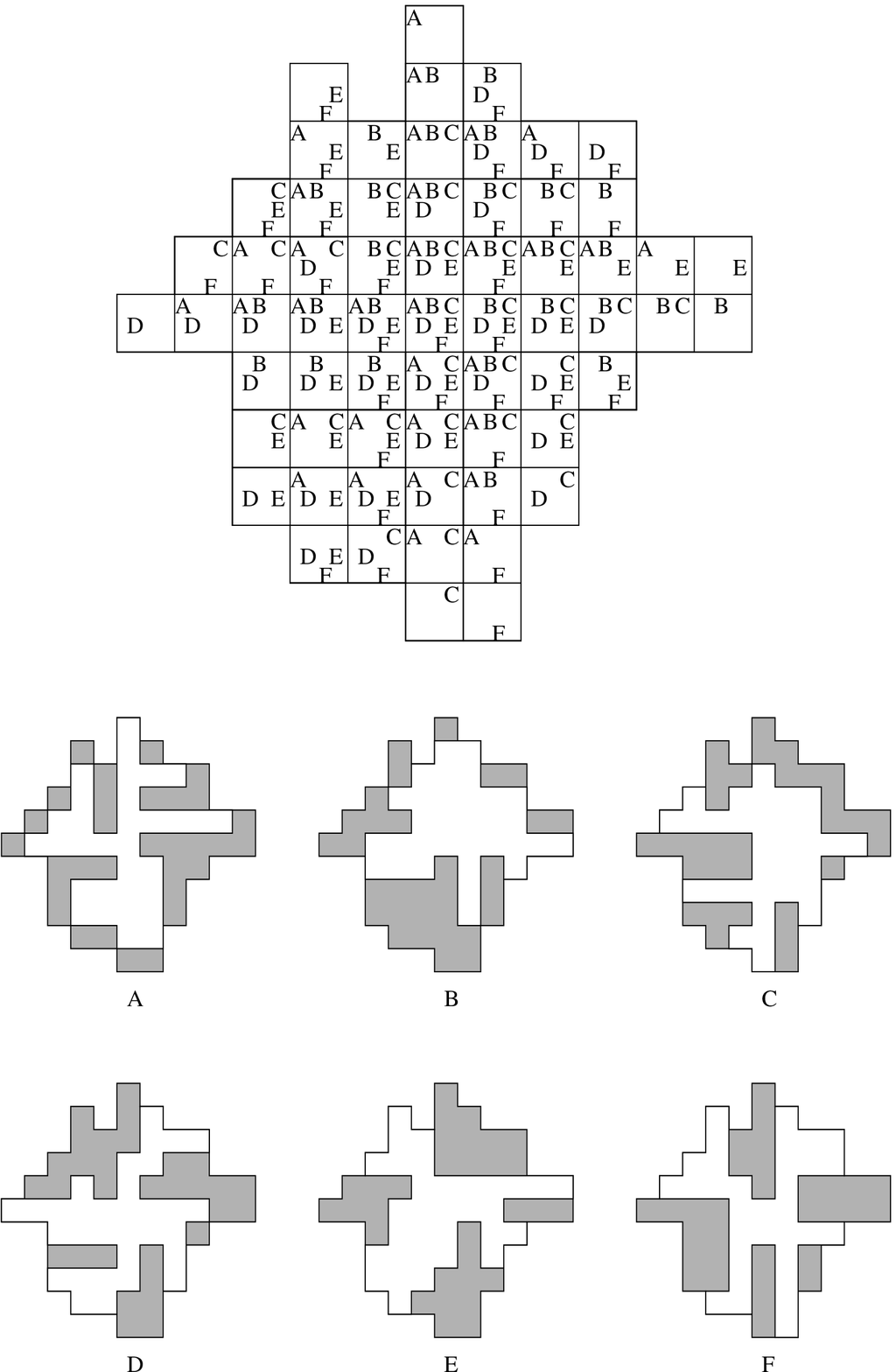}} \caption[]{A
minimum area \polyVenn{6}.} \label{fig:min6}
\end{figure}

\begin{figure}
\centerline{\includegraphics[scale=0.675]{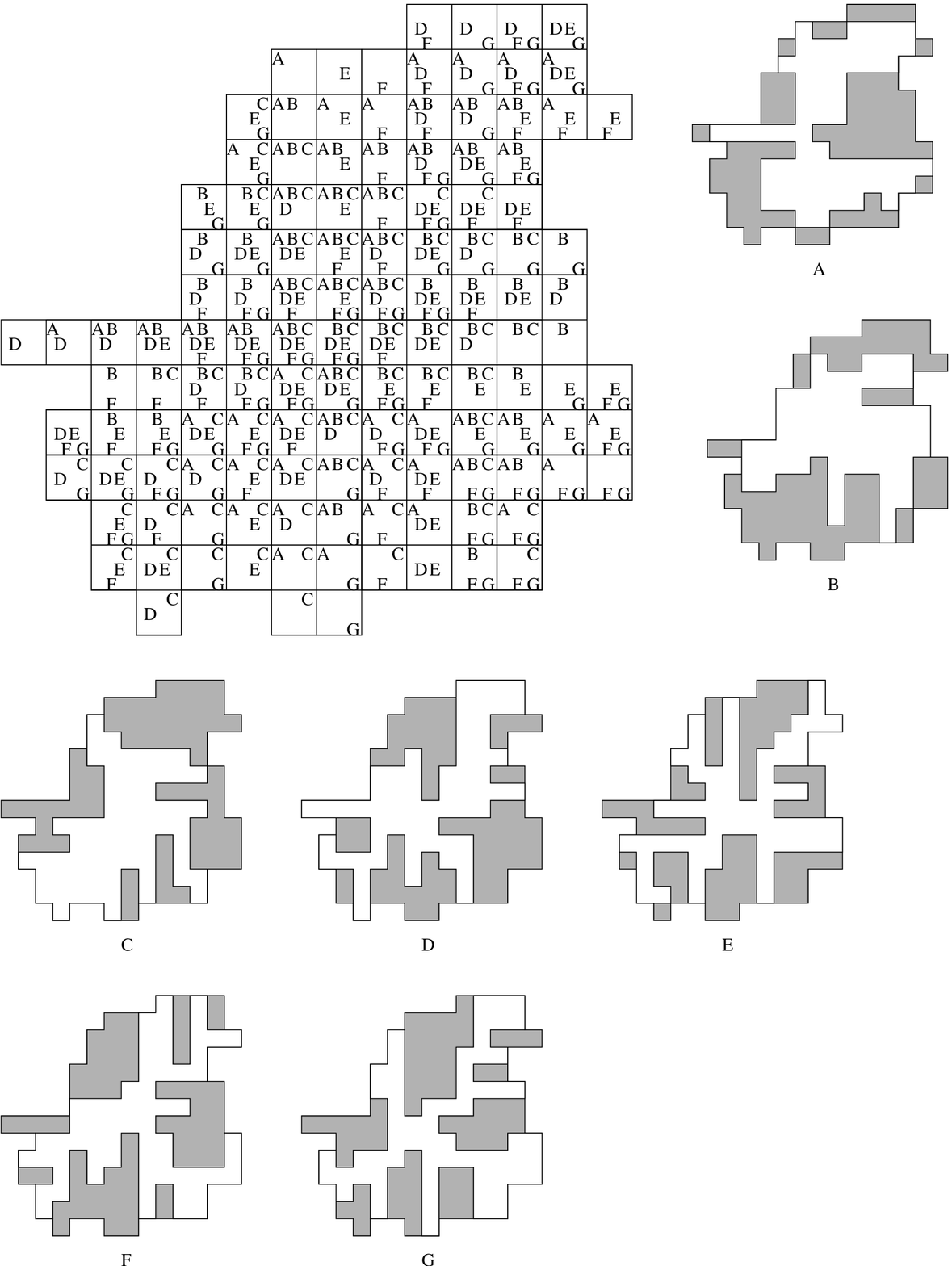}} \caption[]{A
minimum area \polyVenn{7}.} \label{fig:min7}
\end{figure}

By necessity, each curve of a minimum area
  \polyVenn{n} has area $2^{n-1}$.
All the Venn diagrams in Fig. \ref{fig:venn_polyominoes} are
  minimum area congruent \polyVenns{n}.
By trial-and-error, we have also found minimum area
  non-congruent \linebreak[4]
  \polyVenns{n} for $n=6,7$ (see Figs. \ref{fig:min6},\ref{fig:min7}).
It is unknown if minimum area \polyVenns{n}\ exist
  for $n \geq 8$, although we suspect
  there is an upper limit due to the rigid constraints of
  orthogonal grid drawings.

Orthogonal grid drawings of Venn diagrams were first
  studied by Eloff and van Zijl \cite{Eloff};
  they developed a heuristic algorithm based on a
  greedy incremental approach.
An optimization step in the algorithm attempted to
  reduce the overall area of the diagram, but
  there was no upper bound.
In addition, their algorithm produced polyominoes
  with holes, so the resulting diagrams
  would not be considered Venn diagrams in
  the formal sense (because the sets were not
  represented by simple, closed curves).

In the following sections, we present algorithms
  for approximating minimum area \polyVenns{n}.
The first algorithm is trivial and produces \polyVenns{n}\ with
  less than $3/2$ times the minimum area.
The second algorithm improves upon the first by using symmetric
  chain decompositions of the Boolean lattice and produces
  \polyVenns{n}\ whose areas are asymptotically minimum
  (i.e., the ratio of total cells to required cells tends to one
  as $n$ increases).

%

There is another definition of area
  based on the $w \times h$ bounding box that contains an
  \polyVenn{n}; such a box must also have at least one cell to
  represent the empty set.
For example, the \polyVenns{n}\ in Fig.
  \ref{fig:venn_polyominoes} are contained by
  $4 \times 1$,
  $2 \times 5$,
  $5 \times 5$, and
  $7 \times 7$ bounding boxes, respectively.
Since an \polyVenn{n}\ must be comprised of at least
  $2^n-1$ cells, a bounding box must have area at least
  $2^n$.
This leads us to the following definition of a minimum bounding
  box \polyVenn{n}:

\vspace{12pt}
\textsc{Definition} \linebreak \indent A
\textit{minimum bounding box \polyVenn{n}} is an
  orthogonal unit-grid drawing of a Venn diagram that is
  enclosed by a $2^s \times 2^t$ rectangle \linebreak
  where $s+t=n$.
\vspace{12pt}

\begin{figure}
\centerline{\includegraphics[scale=0.8]{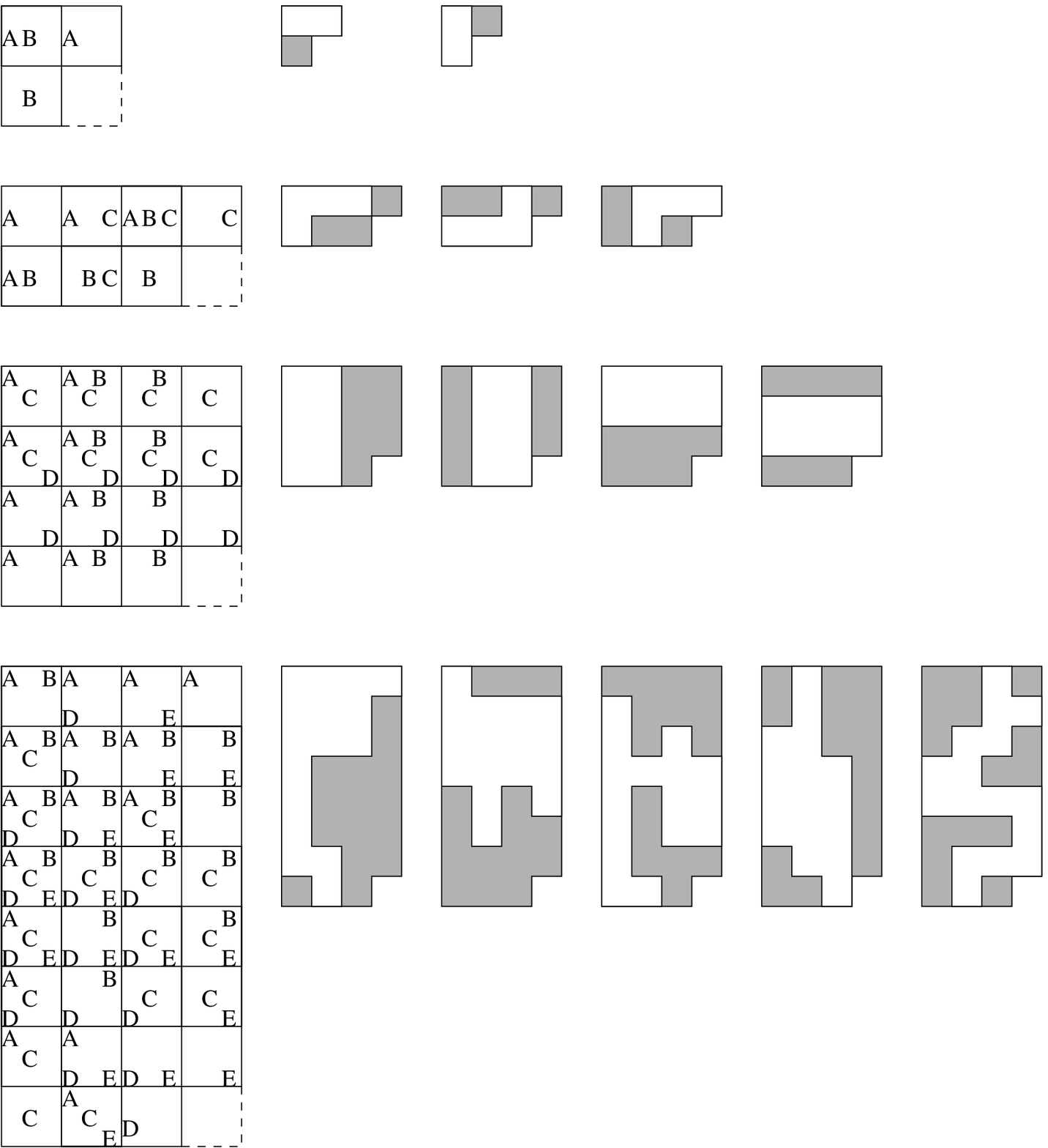}} \caption[
]{Minimum bounding box \polyVenns{n}\ for $2 \le n \le 5$.}
\label{fig:minbbox}
\end{figure}

Of the congruent \polyVenns{n}\ in Fig. \ref{fig:venn_polyominoes},
  only (a) is a minimum bounding box \polyVenn{n}.
Figure \ref{fig:minbbox} shows some examples of minimum bounding
  box non-congruent \polyVenns{n}.

At present, we leave minimum bounding box \polyVenns{n}\ and
  focus the rest of this paper on minimum area \polyVenns{n}.

\subsection*{A $3/2$-APPROX Algorithm}

\begin{figure}
\centerline{\includegraphics[scale=0.875]{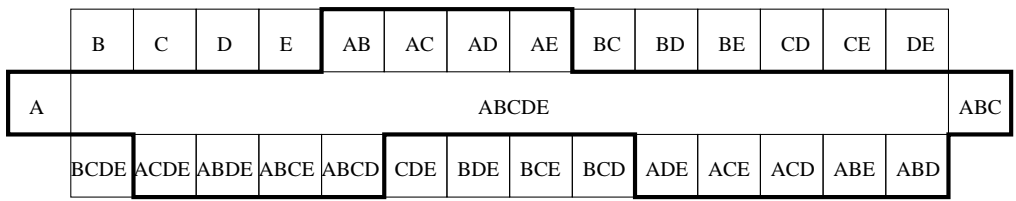}} \caption[
]{A na\"{i}ve approximation for a minimum area \polyVenn{5}; curve
$A$ is highlighted} \label{fig:venn5_naive}
\end{figure}

This algorithm is best explained by way of an example.
Suppose we
  wish to draw a \polyVenn{5}\ with the curves
  $\{A,B,C,D,E\}$.
We begin by drawing a $1 \times 14$ rectangle and labelling it as
  region $ABCDE$;
  in other words, the curves are $1 \times 14$ rectangles stacked
  on top of each other.
We now place $30$ cells around the perimeter of $ABCDE$ and
  uniquely label them with the $30$ remaining non-empty regions;
  the result is shown in Fig. \ref{fig:venn5_naive}.
After adding the perimeter cells, each curve becomes a polyomino
  formed by the original $1 \times 14$ rectangle with ``bumps''
  wherever the curve encloses a perimeter cell.

In the general case, this algorithm will produce an $n$-Venn
  polyomino beginning with a $1 \times (2^{n-1}-2)$ rectangle that has
  a perimeter of $2^n - 2$ (for the $2^n$ regions less the empty and
  full sets).
The resulting diagrams have an area of $2^n + 2^{n-1} - 4$ which
  is less than $3/2$ times the minimum area of $2^n - 1$.

\subsection*{An Asymptotically Optimal Algorithm}

\begin{figure}
\centerline{\includegraphics[scale=1.0]{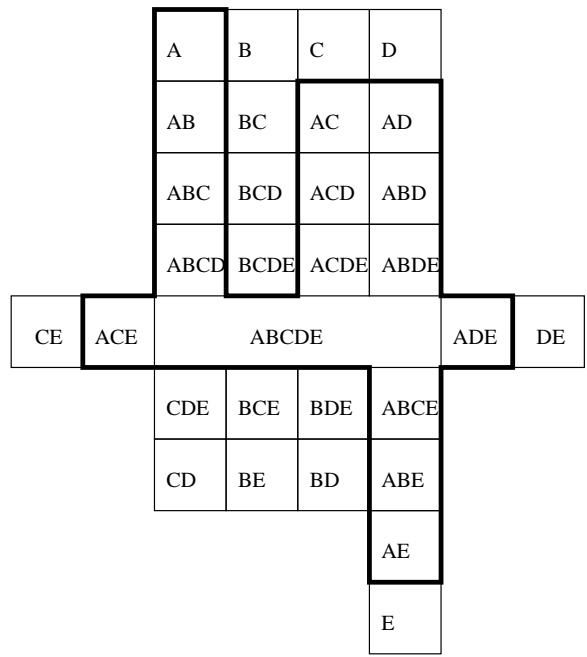}} \caption[
]{An approximation for a minimum area \polyVenn{5}\ using
column-convex polyominoes and symmetric chains; curve $A$ is
highlighted.} \label{fig:venn5_chains}
\end{figure}

The previous algorithm can be significantly improved by noting
  that not all regions need to be placed adjacent to the initial
  rectangle;
  instead, if region $X$ is a subset of region $Y$ then
  $X$ can be placed directly above or below $Y$ (depending on if $Y$
  is above or below the initial rectangle), and the curves will
  remain as polyomino perimeters.
This chaining of regions can continue as long as the subset
  property is maintained.
Figure \ref{fig:venn5_chains} shows an example
  of \polyVenn{5}\ that chains regions as much as possible.
Note also that the resulting polyominoes are column-convex.

When regions are chained, a smaller perimeter is needed for the
  initial rectangle and so the total area of the diagram is reduced.
A smaller area diagram is created by minimizing the number of
  chains, so the question arises as to the best way to decompose the
  regions into chains;
  for this question, we need to use a result from the theory of
  partially ordered sets.

Given a set $S$ with powerset $\powerset{S}$, we define the
  partially ordered set (poset) $\lattice{S}$ with elements
  $\powerset{S}$ ordered by inclusion.
Since $\lattice{S}$ is closed under union, intersection, and
  complement, it is a Boolean lattice.
Figure \ref{fig:poset}(a) shows an example of
  $\lattice{\{A,B,C,D\}}$.

\begin{figure}
\centering \scalebox{0.75}{\input{poset.pstex_t}} \caption[]{(a) A
Hasse diagram of the poset $\lattice{\{A,B,C,D\}}$, (b) one of its
symmetric chain decompositions, and (c) the resulting
\polyVenn{4}.} \label{fig:poset}
\end{figure}

Let $|S|=n$. A symmetric chain decomposition (SCD) of $\lattice{S}$
  is a partition of $S$ into $\binom{n}{\lfloor n/2 \rfloor)}$
  symmetric chains.
Each symmetric chain is a sequence of subsets
  $x_1, x_2, \ldots, x_t$ with the following properties:
\pagebreak
\begin{equation}
x_i \subset x_{i+1} \mbox{ for all } 1 \leq i < t,
\label{eq:subsetprop}
\end{equation}
\begin{equation}
|x_i|=n-|x_{t-i+1}| \mbox{ for all } 1 \leq i \leq \lceil t/2 \rceil.
\end{equation}

Symmetric chain decompositions form an essential ingredient of the
  recent proof of Griggs, Killian and Savage \cite{Griggs} that
  symmetric Venn diagrams exist if and only if the number of curves
  is prime.

Several algorithms exist for decomposing $\lattice{S}$ into
  symmetric chains; we describe two of these algorithms below.
The first, due to de Bruijn, van Ebbenhorst Tengbergen, and
  Kruyswijk \cite{DeBruijn} is called the \emph{Christmas tree pattern}
  by Knuth \cite{Knuth7216}.
It is an inductive construction that creates a set $T_n$ of
  ${n \choose \lfloor n/2 \rfloor}$ chains.
Initially $T_1 = \{ \emptyset \subset \{1\} \}$.  To obtain $T_{n}$ from $T_{n-1}$,
  take each chain $x_1 \subset x_2 \subset \cdots \subset x_t$ in
  $T_n$ and replace it with the two chains
  $x_2 \subset \cdots \subset x_t$ and
  $x_1 \subset x_1 \cup \{n\} \subset x_2 \cup \{n\} \subset \cdots
  \subset x_t \cup \{n\}$ in $T_{n+1}$ if $t > 1$.
If $t=1$ the first chain is empty and is ignored.

A second method, due to Aigner \cite{Aigner}, can be described as
  a greedy lexicographic algorithm.
It is efficient and easy-to-implement, and is the method that we
  used in creating the example diagrams.
Let $m(x,y)$ be the smallest element in a set $x$ that is not in the
  set $y$, where $m(x,y) = -\infty$ if $x \subset y$.
We say that $x$ is \emph{lexicographically smaller} than $y$
  if $m(x,y) < m(y,x)$.
In Aigner's algorithm, the following process is repeated until every
  element of $\lattice{\{1,2,\ldots,n\}}$ is contained in some chain.
For $k = 0,1,2,\ldots,n$, denote by $R(k)$ the set of subsets of $\{1,2,\ldots,n\}$
  size $k$ that are not yet in any chain.
Let $j$ be the smallest value for which $R(j)$ is non-empty and
  let $x$ be the lexicographically smallest set in $R(j)$.
The set $x$ becomes the smallest set in a new chain
  $x = x_1 \subset x_2 \subset \cdots \subset x_t$.
The successive elements of this chain are obtained by taking
  $x_{i+1} \in R(i+1)$ to be the lexicographically smallest set that
  contains $x_{i}$.
It is by no means obvious that this algorithm is correct, but
  indeed it is!

Because of their subset property (\ref{eq:subsetprop}),
  the symmetric chains can be
  directly used to layout the regions of an \polyVenn{n}.
Figure \ref{fig:poset}(b) shows the SCD of $\lattice{\{A,B,C,D\}}$
  that is produced by Aigner's algorithm, and
  Fig. \ref{fig:poset}(c) shows the resulting \polyVenn{4}.
The \polyVenn{5}\ in Fig. \ref{fig:venn5_chains} was also
  produced from Aigner's SCD of $\lattice{\{A,B,C,D,E\}}$.

In the general case, this algorithm will produce an
  \polyVenn{n} beginning with a
  $1 \times (\binom{n}{\lfloor n/2 \rfloor}-2)/2$
  rectangle that has a perimeter of
  $\binom{n}{\lfloor n/2 \rfloor}$.
The resulting diagrams have an area of
  $(\binom{n}{\lfloor n/2 \rfloor}-2)/2 + 2^n - 2$.
The lower bound
  $\binom{2n}{n} < \frac{2^{2n}}{\sqrt{\pi}(n^2+n/2+3/32)^{1/4}}$
  \cite{Grosswald} can be used to show that the algorithm
  produces diagrams whose area is
  $1 + O(1/\sqrt{n})$
  times the minimum area of
  $2^n - 1$;
  therefore, as $n$ increases, the approximation gets
  asymptotically close to optimal.

\subsection*{Open Problems and Final Remarks}

To close the paper, we list some open problems that are inspired by the
  examples in this paper.
With the exception of the congruent \polyVenns{n}, the
  examples in this paper were constructed by hand,
  and it is very likely that relatively na\"ive
  programs will be able to extend them.
Such extension would be interesting, but even more interesting would
  be general results that apply for arbitrary numbers of curves.

\begin{enumerate}
\item Are there congruent \polyVenns{n}\ for $n \geq 6$?
  Figure \ref{fig:venn_polyominoes} shows that they exist for $n = 2,3,4,5$.
\item Is there a \polyVenn{5}\ whose curves are convex polyominoes?
  (The curves in Figure \ref{fig:venn_polyominoes}(d) are not both
  row-convex and column-convex polyominoes.)
\item Are there minimum bounding box \polyVenns{n}\ for $n \geq 6$?
  Figure \ref{fig:minbbox} shows that they exist for $n = 2,3,4,5$.
\item Are there minimum area \polyVenns{n}\ for $n \geq 8$?
  Figure \ref{fig:min7} shows one for $n = 7$.
\item One problem for which we have not attempted solutions
  is the construction of \polyVenns{n}\ that fill a $w \times h$ box,
  where $wh = 2^n-1$.
Of course, a necessary condition is that $2^n-1$ not
  be a Mersenne prime.
For example, is there a \polyVenn{4}\ that fits in a
  $3 \times 5$ rectangle or a \polyVenn{6}\ that fits
  in a $7 \times 9$ or $3 \times 27$ rectangle?
\end{enumerate}

\emph{AUTHORS' NOTE:}
Since submitting the original manuscript of
  this paper, Bette Bultena has discovered a \polyVenn{6} with an
  $8 \times 8$ bounding box (see problem $3$ above).
Figure \ref{fig:min6} has also been used to represent
  the results of experiments in plant genetics \cite{Casimiro}.


\begin{thebibliography}{25}

\bibitem{Aigner} Martin Aigner,
  Lexicographic matching in Boolean algebras,
  \textit{Journal of Combinatorial Theory, Series B}
  \textbf{14} (1973), pp. 187--194.

\bibitem{Battista} Guisseppi Battista, Peter Eades, Roberto Tomassia, and Ionnis Tollis,
  \textit{Graph Drawing: Algorithms for the Visualization of Graphs},
  Prentice-Hall, 1999.

\bibitem{Carroll} Jeremy Carroll,
  6-Venn Triangle Problem,
  \textit{Personal Homepage},
  \verb+http://www.hpl.hp.com/personal/jjc/index.html+.

\bibitem{Casimiro} Sandra Casimiro, Rog\'{e}rio Tenreiro, and
  Ant\'{o}nio A. Monteiro,
  Identification of pathogenesis-related genes in the crucifer downy
  mildew oomycete Hyaloperonospora parasitica by differential display
  analysis of distinct phenotypic interactions with Brassica oleracea,
  \textit{Journal of Microbiological Methods},
  submitted.

\bibitem{DeBruijn} N.G. de Bruijn, C. van Ebbenhorst Tengbergen, and D. Kruyswijk,
  On the set of divisors of a number,
  \textit{Nieuw Archief voor Wiskunde}
  \textbf{2} (1951), pp. 191--193.

\bibitem{Edwards1} Anthony W. F. Edwards,
  Venn diagrams for many sets,
  \textit{New Scientist}
  \textbf{7} (1989), pp. 51--56.

\bibitem{Edwards2} Anthony W. F. Edwards,
  \textit{Cogwheels of the Mind: The Story of Venn Diagrams},
  The Johns Hopkins University Press, 2004.

\bibitem{Eloff} Jacques Eloff and Lynette van Zijl,
  An incremental construction algorithm for Venn diagrams,
  \textit{SAICSIT 2000},
  Cape Town, South Africa, November 2000.

\bibitem{Euler} Leonard Euler,
  \textit{Lettres a une Princesse d'Allemagne sur Divers Sujets de Physique et de Philosophie},
  Imperial Academy of Sciences, Petersburg, 1768--1772.

\bibitem{Golomb} Solomon Golomb,
  \textit{Polyominoes: puzzles, patterns, problems, and packings},
  Princeton University Press, 1994.

\bibitem{Griggs} Jerrold Griggs, Charles E. Killian, and Carla D. Savage,
    Venn diagrams and symmetric chain decompositions in the Boolean lattice,
    \textit{Electronic Journal of Combinatorics},
    \textbf{11} (2004).

\bibitem{Grosswald} E. Grosswald,
  Solution of advanced problem 6019,
  \textit{American Math. Monthly}
  \textbf{84} (January 1977), pp. 63--65.

\bibitem{Grunbaum0} Branko Gr\"{u}nbaum,
  Venn diagrams and Independent Families of Sets,
  \textit{Mathematics Magazine}
  \textbf{48} (Jan-Feb 1975), pp. 12--23.

\bibitem{Grunbaum1} Branko Gr\"{u}nbaum,
  Venn diagrams I,
  \textit{Geombinatorics} (Vol. I)
  \textbf{4} (1992), pp. 5--12.

\bibitem{Grunbaum2} Branko Gr\"{u}nbaum,
  Venn diagrams II,
  \textit{Geombinatorics} (Vol. II)
  \textbf{2} (1992), pp. 25--32.

\bibitem{Grunbaum3} Branko Gr\"{u}nbaum and Peter Winkler,
  A Venn diagram of 5 triangles,
  \textit{Mathematics Magazine}
  \textbf{55} (1982), pg. 311.

\bibitem{James} Ioan James,
  \textit{Remarkable Mathematicians: From Euler to von Neumann},
  Cambridge University Press, 2002.

\bibitem{Jensen} Iwan Jensen and Anthony J. Guttmann,
  Statistics of lattice animals (polyominoes) and polygons,
  \textit{Journal of Physics, Series A}
  \textbf{33} (2000), pp. L257--L263.

\bibitem{Klarner} David A. Klarner,
  Some results concerning polyominoes,
  \textit{Fibonacci Quarterly}
  \textbf{3} (1965), pp. 9--20.

\bibitem{Knuth7216} Donald E. Knuth,
  \textit{The Art of Computer Programming},
  pre-fascicle 4A (a draft of Section 7.2.1.6: generating all trees),
  Addison-Wesley, 2004.

\bibitem{Mertens} Stephan Mertens,
  Counting lattice animals: a parallel attack,
  \textit{Journal of Statistical Physics}
  \textbf{66} (1992), pp. 669--678.

\bibitem{Polya} George P\'{o}lya,
  On the number of certain lattice polygons,
  \textit{Journal of Combinatorial Theory}
  \textbf{6} (1969), pp. 102--105.

\bibitem{Redelmeier} D.H. Redelmeier,
  Counting polyominoes: yet another attack,
  \textit{Discrete Mathematics}
  \textbf{36} (1981), pp. 191--203.

\bibitem{Ruskey} Frank Ruskey,
  A survey of Venn diagrams,
  \textit{Electronic Journal of Combinatorics}
  \textbf{4} (1997), DS \#5 (updated 2001).

\bibitem{Sloane} Neil J.A. Sloane,
  \textit{The On-Line Encyclopedia of Integer Sequences},
  \verb+http://www.research.att.com/~njas/sequences/+.

\bibitem{Thompson} Mark Thompson,
  Venn polyominoes,
  \textit{Math Recreations},
  \verb+http://www.flash.net/~markthom/html/venn_polyominos.html+.

\bibitem{Venn} John Venn,
  On the diagrammatic and mechanical representation of propositions and reasonsings,
  \textit{The London, Edinburgh, and Dublin Philosophical Magazine and Journal of Science}
  \textbf{9} (1880), pp. 1--18.

\end{thebibliography}
\end{document}